\documentclass{amsart}[12pt]
\usepackage{amsmath,mathrsfs,amssymb}
\usepackage{fullpage}

\newtheorem{lem}{Lemma}

\newtheorem*{rem}{Remark}
\title{The description of the free lognormal distribution revisited} 
\author[N. Demni]{Nizar Demni}
\address{IRMAR, Universit\'e de Rennes 1\\ Campus de
Beaulieu\\ 35042 Rennes cedex\\ France}
\email{nizar.demni@univ-rennes1.fr}

\author[T. Hamdi]{Tarek Hamdi}
\address{IPEST, Universit\'e de Carthage\\ 2078 La Marsa \\ Tunisie}
\email{tarek.hamdi@ipest.rnu.tn}
\begin{document}
\maketitle
\begin{abstract}
We adapt the approach used in \cite{DH} and retrieve the description given in \cite{Biane1} of the free lognormal distribution. 
\end{abstract}

\section{Introduction}
Let $t \in \mathbb{R}$ and define the sequence
\begin{equation}\label{M}
m_n(t) := \frac{e^{-nt/2}}{n} \sum_{k=0}^{n-1}\frac{(-tn)^k}{k!}\binom{n}{k+1}, \quad n \geq 1.  
\end{equation}
It is not difficult to see that $m_n(t)$ can be written through the $(n-1)$-th Laguerre polynomial as follows (\cite{Szego}): 
\begin{equation*}
m_n(t) = \frac{e^{-nt/2}}{n}L_{n-1}^{(1)}(nt). 
\end{equation*}
The analysis of the large-$n$ behaviour of this kind of polynomials-for which the argument and the degree are dependent- appeared for instance in \cite{Szego}, \cite{LT} and \cite{Gro-Mat} and relies on the saddle point method. However, what is far from being clear and true is the fact that $(m_n(t))_{n \geq 0}$ with $m_0(t) := 1$ is the moment sequence of a compactly-supported probability distribution $\nu_t$ (\cite{Biane}). More amazing is that $\nu_t$ is supported in the unit circle when $t > 0$ and in an interval in the positive real line otherwise (\cite{Biane1}). Note that a similar phenomenon happens when the index and the argument of the Laguerre polynomial are independent. Indeed, it is already known that 
\begin{equation*}
n! L_n^{(\alpha)}(x), \,\, n \geq 0,
\end{equation*}
is the moment sequence of a positive measure supported in the whole positive half-line provided that $x < 0$ (\cite{IS}, 523-524). On the other hand, if $x > 0$ and $\alpha > -1/2$, then the Laguerre polynomial admits a Laplace-type integral representation over the whole upper half-plane which follows from the well-known Laplace-type integral representation of Jacobi polynomials due to Koornwinder (see e.g. \cite{AAR} p.475 and formula 18.10.6 in \cite{DLMF}). In particular, $n!L_n^{(\alpha)}$ is the $n$-th `complex' moment of a positive measure. 

Coming back to $(m_n(t))_n$, the fact that it is the moment sequence of a probability distribution follows mainly from the large-size asymptotics of observables of the Brownian motion on Lie groups. For instance, let $(Y_t)_{t \geq 0}$ be a Brownian motion on the unitary group $U(d)$ (\cite{Liao}). Then, for any $t \geq 0$ and any $n \in \mathbb{N}$,
\begin{equation*}
m_n(t) = \lim_{d \rightarrow \infty} \frac{1}{d}\mathbb{E}(\textrm{tr}(Y_{t/d}^n)).  
\end{equation*}
Moreover, $Y_t^{\star}$ and $Y_t$ have the same distribution so that $m_{-n}(t) = m_n(t), n \geq 0$. This result, proved independently in \cite{Biane}, \cite{Rains} and \cite{Xu}, was further strengthened in \cite{Levy} where a full expansion of the averaged power sums of $Y_t$ is obtained and involves paths in the Cayley graph of the symmetric group. As to the density and the support of $\nu_t, t \geq 0$, they were described in \cite{Biane1} from the behaviour of the compositional inverse of the Herglotz transform of $\nu_t$ in some Jordan domain. In particular, the support exhibits a phase transition at $t = 4$ which agrees with the change of the decay of $m_n(t)$ for large $n$. 
For times $t \in (0,4)$, the same description as well as the phase transition were found in \cite{DH} relying on a suitable Cauchy-type integral representation of the Laguerre polynomials: the support is built from the intersection of two curves which disconnect exactly at $t=4$. 

When $t \leq 0$, the moments $(m_n(t))_n$ turn out to be related to the radial part of the right-invariant Brownian motion on the linear complex group $GL(d,\mathbb{C})$ (\cite{Liao}, \cite{NRW}). More precisely, let $(Z_s)_{s \geq 0}$ be such a process, then it was conjectured in \cite{Biane} and proved later in \cite{Ceb} that
\begin{equation}\label{LR2}
m_n(t) = \lim_{d \rightarrow \infty} \frac{1}{d}\mathbb{E}(\textrm{tr}[(Z_{-t/d}^{\star}Z_{-t/d})^n)]),\, n \geq 0.
\end{equation}
Moreover, $\nu_t$ is compactly-supported in the positive half-line and can be realized as the free additive convolution of a semi-circular distribution of variance $-t$ and a uniform distribution on $[t/2,-t/2]$ (\cite{Ho}). It is also called the free log-normal distribution since it appears as the central limit of the products of positive free random variables (\cite{Ho}, \cite{SY}). 

In this paper, we adapt the approach used in \cite{DH} for the description of $\nu_t, t \geq 0$ to the free lognormal distribution $\nu_t, t \leq 0$. More precisely, using a Cauchy-type integral representation of Laguerre polynomials, we prove that for any $t < 0$, there exists a Jordan curve $\gamma_t$ around the origin whose image under the (meromorphic) integrand lies in the positive half-line. The description of the free lognormal distribution follows then from two variable changes performed on two parts of the curve where the integrand is shown to be monotone. 

\section{Description of $\nu_t, t \geq 0$ revisited} 
\subsection{The curve $\gamma_t$} 
Since   
\begin{equation*}
\int_{\gamma} e^{-ntz}\frac{dz}{z^k} = \frac{(-nt)^{k-1}}{(k-1)!}, \, 1 \leq k \leq n, 
\end{equation*}
and vansihes for $k=0$ then
\begin{equation*}
m_n(t) = \frac{e^{-nt/2}}{2i\pi n}\int_{\gamma} e^{-ntz}\left(1+\frac{1}{z}\right)^n dz.
\end{equation*}
This is nothing else but one of the numerous integral representations of Laguerre polynomials (\cite{Szego}) and it was already used in \cite{Gro-Mat} in order to estimate the growth of $m_n(t), t > 0$, for large $n$. Now, we already know that $\nu_t$ is compactly supported in the positive half-line. Hence, it is natural to seek for any fixed $t < 0$ a Jordan curve around the origin $\gamma_t$ such that 
\begin{equation*}
g_t(z) := e^{-t(z+(1/2))}\left(1+\frac{1}{z}\right), \quad z \in \gamma_t.
\in \mathbb{R}
\end{equation*}  
In this respect, we easily prove the following:
\begin{lem}
For any $t < 0$, there exists a Jordan curve $\gamma_t$ around the origin such that $g_t(z) \geq 0, z \in \gamma_t$. 
\end{lem}
\begin{proof}
Write $z = x+iy$, then 
\begin{equation}\label{Const1}
g_t(z) \in \mathbb{R} \quad \Leftrightarrow \quad y\cos(ty) + (x^2+y^2+x)\sin(ty) = 0.
\end{equation} 
Since \eqref{Const1} is satisfied by any real number, then we rather focus on 
\begin{equation}\label{Const2}
-y^2-y\cot(ty) = x^2+x.
\end{equation}
But $\cot(\pm \pi) = \infty$ so that \eqref{Const2} is satisfied as soon as one finds a positive number $0 < y_t < -\pi/t$ such that
\begin{equation*}
\frac{1}{4} - y^2 - y\cot(ty) \geq 0 
\end{equation*}
for all $|y| \leq y_t$. Hence, the curve $\gamma_t$ we are looking for would be defined by  
\begin{equation*}
x_t^{\pm}(y) = -\frac{1}{2} \pm \sqrt{\frac{1}{4} - y^2 - y\cot(ty)}
\end{equation*}
and meet the real line at the two points
\begin{equation*}
x_t^{\pm}(0) = -\frac{1}{2} \pm \sqrt{\frac{1}{4} - \frac{1}{t}}. 
\end{equation*}  
Now, set
\begin{equation*}
f_t(y) := 2y\sin(ty) + \cos(ty), \, y \in (\pi/t, -\pi/t),
\end{equation*}
then for any $y \neq 0$
\begin{equation*}
\frac{1}{4} - y^2 - y\cot(ty) \geq 0 \quad \Leftrightarrow \quad |f_t(y)| \leq 1.
\end{equation*}
Let $y \in [0,-\pi/t)$, then the inequality $\sin u \geq u, u \leq 0$ shows that
\begin{equation*}
y \mapsto (ty)\cot(ty)
\end{equation*}
is decreasing. But 
\begin{equation*}
f_t'(y) = 2\sin(ty)\left[ty\cot(ty) - \frac{t-2}{2}\right] 
\end{equation*}
therefore there exists $a_t \in [-\pi/(2t),-\pi/t[$ such that $f_t'(a_t) = 0$ whence  \\

\centerline{$
\begin{array}{|c|ccccr|}
\hline
y     & 0   & & a_t  &    & -\pi/t  \\
\hline
f_t'(y) &  - &    & 0  &  & + \\
\hline &&  &&   & \\       
f_t(y) &1  &\searrow  & f_t(a_t) &\nearrow  & -1 \\   
&  & & & &  \\         
\hline
\end{array} \quad .
$
}
\vspace{0.2cm}
Besides, the very definition of $a_t$ implies that 
\begin{equation*}
f_t(a_t) = \frac{\sin(ta_t)}{ta_t}[2ta_t^2 + \frac{t}{2}-1] \leq 2ta_t^2 + \frac{t}{2} -1 < -1.
\end{equation*}
Consequently, there exists $0 < y_t < a_t$ such that $|f_t(y)| \leq 1$ for all $y \in [-y_t,y_t]$ and as such, the curve $\gamma_t$ exists. 
Finally, for any $z \in \gamma_t$,
\begin{align*}
g_t(z) = \Re(g_t(z)) &= \frac{e^{-t(x+(1/2))}}{x^2+y^2} \left[(x^2+y^2+x)\cos(ty) -y \sin(ty)\right]
\\& =   y\frac{e^{-t(x+(1/2))}}{x^2+y^2} \left[\cot(-ty)\cos(ty) + \sin(-ty)\right]
\\& =  \frac{y}{\sin(-ty)} \frac{e^{-t(x+(1/2))}}{x^2+y^2} 
\\& =  \frac{y}{\sin(ty)} \frac{e^{-t(x+(1/2))}}{y\cot(ty)+x} 
\end{align*}
which is obviously positive. 
\end{proof}
\begin{rem}
The curve $\gamma_t$ is clearly symmetric with respect to the real axis and the vertical line $x=-1/2$. Moreover, the transformations $y \rightarrow y/2, x \rightarrow (x -1)/2$ and $t \rightarrow -4t$ map it onto the boundary of the region $\Omega_t$ described p.271 in \cite{Biane1}: 
 \begin{equation*}
2y\cot(2ty) = x^2+y^2 -1.
\end{equation*} 
\end{rem}
\subsection{Monotonicity of $y \mapsto g_t(x_t^{\pm}(y),y)$}
Recall from the end of the proof of the previous lemma that 
\begin{equation*}
g_t(x,y) =  \frac{y}{\sin(ty)} \frac{e^{-t(x+1/2)}}{y\cot(ty)+x} 
\end{equation*}
for any $z = x+iy \in \gamma_t$. Set
\begin{equation*}
v_t(y) := \sqrt{\frac{1}{4} - y^2 - y\cot(ty)}, \quad y \in [0,y_t],
\end{equation*}
and recall the notation $x_t^{\pm}(y) = (1/2) \pm v_t(y)$. Then 
\begin{equation*}
g_t(x_t^{\pm}(y),y) =  \frac{y}{\sin(ty)} \frac{e^{\mp tv_t(y)}}{y\cot(ty) - (1/2) \pm v_t(y)}. 
\end{equation*}
In this paragraph, we shall prove that 
\begin{lem}
For any $t < 0$, the map $y \mapsto g_t(x_t^-(y),y)$ (respectively $y \mapsto g_t(x_t^+(y),y)$) is increasing (respectively decreasing) on $(0,y_t)$. 
\end{lem}
\begin{proof}
Firstly, note that 
\begin{equation*}
g_t(x_t^-(y),y)g_t(x_t^+(y),y) = 1,
\end{equation*}
therefore it suffices to prove that $y \mapsto g_t(x_t^-(y),y)$ is increasing. Now, for any $y \in (0,y_t)$,
\begin{equation*}
\partial_yg_t[x_t^{-}(y),y] = \partial_y x_t^{-}(y)\partial_x(g_t)[x_t^{-}(y),y] + \partial_y(g_t)[x_t^{-}(y),y]
\end{equation*} 
where 
\begin{eqnarray*}
\partial_x(g_t)(x,y) &=& -\frac{ye^{-t(x+1/2)}}{(x+y\cot(ty))^2\sin(ty)} [t(x+y\cot(ty)) + 1] \\
\partial_y(k_t)(x,y) &=& \frac{e^{-t(x+1/2)}}{(x+y\cot(ty))^2\sin(ty)} [(ty^2+x) - txy\cot(ty)]\\ 
\partial_y(x_t^{-})(y) &=&  \frac{4y\sin^2(ty) + \sin(2ty) - 2ty}{4\sin^2(ty)v_t(y)}.
\end{eqnarray*}
Observe that the analysis performed in the previous paragraph shows that $v_t(y) = 0, y \in [0,y_t]$ if and only if $y = y_t$. Besides, the inequality $\sin u \geq u, u \leq 0$ shows that $\partial_y(x_t^{-})(y) > 0$ so that $y \mapsto x_t^-(y)$ is increasing from $(0,y_t)$ onto $(x_t^-(0), -1/2)$. Secondly, the equality
\begin{equation*}
x + y\cot(ty) = -(y^2+x^2)
\end{equation*}
along the curve $\gamma_t$ shows that the derivative $\partial_x(g_t)[x_t^{-}(y),y]$ is positive on $(0,y_t)$. Finally, the positivity of $\partial_y(k_t)(x_t^-(y),y)$  on $(0,y_t)$ follows from the inequalities 
\begin{equation*}
x_t^-(y) < 0 \quad \textrm{and} \quad 1-ty\cot(ty) > 0.
\end{equation*} 
The lemma is proved.
\end{proof}

\subsection{Description of $\nu_t$} 
Since $\gamma_t$ is piecewise smooth and closed, then an integration by parts yields 
\begin{align*}
m_n(t) & =  \frac{1}{2i\pi n}\int_{\gamma_t} [g_t(z)]^n dz = - \frac{1}{2i\pi}\int_{\gamma_t} [g_t(z)]^n z\frac{\partial_zg_t(z)}{g_t(z)}dz.
\end{align*}
Now split $\gamma_t = \gamma_t^+ \cup \gamma_t^-$ where 
\begin{equation*}
\gamma_t^{\pm} := \{z_t^{\pm}(y) := x_t^{\pm}(y)+ iy, \, |y| \leq y_t\}.
\end{equation*}
Hence
\begin{align*}
\int_{\gamma_t} [g_t(z)]^n z\frac{\partial_zg_t(z)}{g_t(z)}dz &= \int_{-y_t}^{y_t}[g_t(z_t^+(y))]^n z_t^+(y)\frac{\partial_y[g_t(z_t^+)](y)}{g_t(z_t^+(y))}dy
\\&  - \int_{-y_t}^{y_t}[g_t(z_t^-(y))]^n z_t^-(y)\frac{\partial_y[g_t(z_t^-)](y)}{g_t(z_t^-(y))}dy.
\end{align*}
But $z_t^+(-y) = \overline{z_t^+(y)}$ ans $g_t(z_t^+(y)) \in \mathbb{R}$ therefore $g_t(z_t^+(-y))= g_t(z_t^+(y))$ which in turn entails 
\begin{align}\label{E1}
 \int_{-y_t}^{y_t}[g_t(z_t^+(y))]^n z_t^+(y)\frac{\partial_y[g_t(z_t^+)](y)}{g_t(z_t^+(y))}dy & =  \int_{-y_t}^{y_t}[g_t(z_t^+(y))]^n z_t^+(y)\partial_y[\log g_t(z_t^+)](y)dy \nonumber
\\& = \int_0^{y_t}[g_t(z_t^+(y))]^n z_t^+(y)\partial_y[\log g_t(z_t^+)](y)dy  \nonumber
\\& -  \int_0^{y_t}[g_t(z_t^+(y))]^n \overline{z_t^+(y)}\partial_y[\log g_t(z_t^+)](y)dy   \nonumber
\\& = 2i \int_0^{y_t}[g_t(z_t^+(y))]^n \Im[z_t^+(y)]\partial_y[\log g_t(z_t^+)](y)dy. 
 \end{align}
Similarly 
\begin{align}\label{E2}
\int_{-y_t}^{y_t}[g_t(z_t^-(y))]^n z_t^-(y)\frac{\partial_y[g_t(z_t^-)](y)}{g_t(z_t^-(y))}dy = 2i \int_0^{y_t}[g_t(z_t^-(y))]^n \Im[z_t^-(y)]\partial_y[\log g_t(z_t^-)](y)dy.
\end{align}
Gathering \eqref{E1} and \eqref{E2}, we get
\begin{align*}
\frac{1}{2i\pi} \int_{\gamma_t^{\pm}}  [g_t(z)]^n z\frac{\partial_zg_t(z)}{g_t(z)}dz &= \frac{1}{\pi}\left\{\int_0^{y_t}[g_t(z_t^+(y))]^n \Im[z_t^+(y)]\partial_y[\log g_t(z_t^+)](y)dy\right. 
\\& - \left. \int_0^{y_t}[g_t(z_t^-(y))]^n \Im[z_t^-(y)]\partial_y[\log g_t(z_t^-)](y)dy\right\}.
\end{align*} 
Finally, the monotonicity of $y \mapsto g_t(z_t^{\pm}(y))$ entails   
\begin{equation*}
m_n(t) = \frac{1}{\pi} \int_{g_t(x_t^-(0))}^{g_t(x_t^+(0))} x^n \Im[g_t^{-1}(x)] \frac{dx}{x}
\end{equation*}
where 
\begin{equation*}
g_t(x_t^{\pm}(0)) = \left[1-\frac{t}{2} \pm \sqrt{\frac{t^2}{4} - t}\right] \exp\{\pm\sqrt{\frac{t^2}{4} - t}\}.
\end{equation*}
Using the time change $t \rightarrow -4t$, the interval $[g_t(x_t^{-}(0)), g_t(x_t^{+}(0))]$ transforms into the support of $\nu_t$ written in \cite{Biane}, Proposition 11.

\end{document}